\documentclass[11pt,a4paper]{article}
\usepackage[T1]{fontenc}
\usepackage{amsmath}
\usepackage{amssymb}

\begin{document}

\centerline{{\bf AFFINE HOMOGENEOUS STRICTLY PSEUDOCONVEX}}
\centerline{{\bf HYPERSURFACES OF THE TYPE $ (1/2,0) $ IN $
\mathbb{C}^3 $}}

\

\centerline{\bf A.V. Atanov$^{1}$, A.V. Loboda$^{1}$, A.V.
Shipovskaya$^{2}$} \centerline{${}^{1}$ Voronezh State University}
\centerline{${}^{2}$ Voronezh State University of Architecture and
Civil Engineering}

\centerline{\textit{Partially supported by RFBR grant 11-01-00495-a
}}

\

\begin{abstract}
In this paper we denote a type of affine homogeneous real
hypersurface of $\mathbb{C}^3$ and
present a classification of homogeneous
surfaces of the type $(1/2, 0) $.
The result was obtained by reducing the classification problem
mentioned above to the problem of solving a system of
nonlinear (quadratic) equations. Solutions of such system represent all the
Lie algebras corresponding to the homogeneous surfaces under consideration.
Some details of integration procedure are discussed for obtained algebras.
\end{abstract}

\

\centerline{\textbf{1. Introduction}}

\

   The description problem for real submanifolds of complex spaces
(CR-manifolds) having "reach" \ symmetry groups, and in particular
homogeneous manifolds, has been widely studied in recent years (see,
for example [1]-[10]).

The first classification result for \emph{holomorphically}
homogeneous real hypersurfaces of 2-dimensional complex spaces (3-dimensional
CR-manifolds) was given by Cartan in [1].
Beloshapka and Kossovskiy completely described 4-dimensional
homogeneous CR-manifolds in [10].

   However, so far only some partial classification results have been obtained
for the next 5-dimensional case of homogeneous CR-manifolds.

   Note that any analytic CR-manifold can be (locally) considered as the one
embedded in some complex space $ \mathbb{C}^n $ of appropriate dimension.
As another obvious remark we
emphasize that most of the mentioned classification results for
holomorphic homogeneity are closely related to the affine
homogeneity property. This remark explains the exceptional interest
to the affine homogeneity in the holomorphic geometry context.

   The
classification of \emph{affine} homogeneous real hypersurfaces of $ \Bbb C^2 $
was obtained in [11].
In this paper we consider affine homogeneity
property for the 5-dimensional real hypersurfaces of 3-dimensional
complex space.

   A real hypersurface $ M $ of $\mathbb{C}^n $
(or $\mathbb{R}^n $) is called {\it affine homogeneous near
a point $ P \in M $} if some Lie subgroup $ G(M) $
of affine group $ Aff(n,\mathbb{C}) $ (or $ Aff(n,\mathbb{R}) $)
acts transitively on $ M $ in a neighborhood of P.

    There are several approaches to the homogeneity problem
for real hypersurfaces in holomorphic as well as in the affine cases.
One of them is the
algebraic approach which is based on the study of canonical
equations of the real analytic hypersurfaces. Below we discuss only
strictly pseudo-convex (SPC) hypersurfaces (one can find the
definition of the strictly pseudo-convexity in [17]).

It was proved
in [13] that in a neighborhood of any point an equation of
real analytic hypersurface $ M \in \Bbb C^3 $
can be written after some affine transformation in the form:

$$\begin{array}{c}
  v = |z_1|^2 + |z_2|^2 + \varepsilon_1 (z_1^2 + \bar z_1^2) +
                           \varepsilon_2 (z_2^2 + \bar z_2^2) + \\
    + \sum\limits_{k+l+2m \ge 3} F_{klm} (z,\bar z) u^m.
\end{array} \eqno (1.1)$$

   Here $ z_1, z_2, w $ are complex coordinates in $ \mathbb{C}^3 $,
$ u = Re \,w, v = Im \,w $;

$F_{klm} $ is a polynomial of degree $ k $ in variable $ z = (z_1,
z_2) $, degree $ l $ in variable $ \bar z = (\bar z_1, \bar z_2)$,
and degree $\ m $ in variable $ u $.

   Unordered pair $(\varepsilon_{1}, \varepsilon_{2})$ of real
nonnegative coefficients from (1.1) is the affine invariant of the
hypersurface at a fixed point. In the affine homogeneous case this
pair will be called a \emph{type} of the homogeneous hypersurface.

   Complete classification of affine homogeneous hypersurfaces of the
type $(1/2,1/2)$ in $\mathbb{C}^3$ was presented in [14] - [16]. The
tubular hypersurfaces (tubes) $ T = \Gamma + i \mathbb{R}^3 $ over
strictly convex affine homogeneous hypersurfaces $ \Gamma \subset
\mathbb{R}^3 $ are the main examples of homogeneous CR-manifolds in
this case. However, it was shown in [14] - [16] that a lot of
homogeneous hypersurfaces of $(1/2, 1/2)$-type cannot be affinely
reduced to the tubes.

   In addition, recall that a big family
of homogeneous hypersurfaces of the type $ (\varepsilon, 0) $ where $ 0<
\varepsilon \ne 1/2 $ was presented in [13].

  In this paper we consider affine homogeneous hypersurfaces of the
type $(1/2, 0) $. Following the scheme introduced in [14], we reduce
the classification problem for affine homogeneous hypersurfaces to
the problem of solving a system of nonlinear equations and construct
a list of affine homogeneous hypersurfaces of the type $(1/2, 0) $.
This list is relatively long and (as the authors hope) complete.
However, it is clear that any classification result containing many
objects needs to be checked by multiple verifications.

\

\centerline{\textbf{2. The main theorem}}

\

{\bf Theorem.} {\it The following hypersurfaces of the type (1/2,0)
are affine homogeneous manifolds in $ \mathbb{C}^3 $ :
$$
   v =  2 x_1^2 + |z_2|^2;
\eqno (2.1)
$$
$$
   v = \exp (x_1) + |z_2|^2,
\eqno (2.2)
$$
$$
   v = - \ln(1+x_1) + |z_2|^2 \ (x_1 > -1),
\eqno (2.3)
$$
$$
    v = \pm (1+x_1)^{\alpha} + |z_2|^2 \ (x_1 > -1),
        \alpha \in \mathbb{R} \setminus \{0,1,2 \},
\eqno (2.4)
$$
$$
     v = (1+x_1)\ln(1+x_1) + |z_2|^2   \ (x_1 > -1),
\eqno (2.5)
$$
$$
    v^2 = |z_1|^2 + |z_2|^2 \ ( v \ne 0),
\eqno (2.6)
$$
$$
      v = \frac{x_1^2}{1-x_2} + |z_2|^2 \ (x_2 \ne 1),
\eqno (2.7)
$$
$$
    v = |z_1|e^{B \arg{z_1}} + |z_2|^2 \ (z_1 \ne 0), \ B \in \mathbb{R};
\eqno (2.8)
$$
$$
     v = x_1^{1-\alpha} |z_2|^{2\alpha} \ (x_1 \cdot |z_2| > 0), \
                    \alpha\in \mathbb{R} \setminus\{0,1\};
\eqno (2.9)
$$
$$
    Re  (\bar z_1 w) =  (Re (z_1 \bar z_2))^{\alpha} \
    ( Re (z_1 \bar z_2) > 0), \
             \alpha \in (-\infty,0).
\eqno (2.10)
$$
}

{\bf Remark.} The $ \pm $ sign in equation (2.4) depends on the
parameter $ \alpha $. The required sign should ensure the
positiveness of a coefficient under the term $ x_1^2 $ in a Taylor
expansion of the right-hand side of (2.4).

\

  The main theorem has a certain connection with the problem of
(local) classification of holomorphically homogeneous hypersurfaces
of $\mathbb{C}^3$.

   For example, hypersurfaces (2.1)-(2.5) are well-known tubular
hypersurfaces (tubes); hypersurface (2.8) is locally holomorphically
equivalent to a sphere in $\mathbb{C}^3$; hypersurfaces (2.6), (2.9)
and (2.10) are holomorphically equivalent to a non-spherical tubes
$$
    v = \ln x_1 + A \ln x_2 \  (A \ne 1)
$$
with 7-dimensional groups of holomorphic transformations.

  A cubic hypersurface (2.7) is the most interesting among the listed manifolds.
It is an element of 1-parameter family of affine homogeneous
hypersurfaces
$$
      v = \frac{x_1^2}{x_2} + ( A x_2^2 + y_2^2) \ (A \in \mathbb{R}).
$$

  Note that any such surface has a type $ (1/2, (A-1)/2(A+1)) $ if $ A \ne -1 $.
So the simple generalization of the example (2.7) gives a big family
of different affine homogeneous hypersurfaces.

\

  All the hypersurfaces listed in the main theorem were obtained by
constructing and integrating a large number of 5-dimensional Lie
algebras related to homogeneous hypersurfaces of the type (1/2,0).

  Any such algebra $ g(M) $
can be considered as the algebra of affine vector fields tangent to
homogeneous hypersurface $ M $.
 The values of these fields
at a point $P \in M$ cover all the tangent plane $ T_P M $.

   It can be shown that every algebra $g(M)$ has a
matrix representation, so we associate with any vector field
$$
{\begin{array}{*{20}c}
 {Z=(a_{1} z_1 + a_{2} z_2 + a_{3} w + p)\frac{\partial }{\partial z_1 }+}
\hfill \\
\qquad
 {+(b_{1} z_1 + b_{2} z_2 + a_{3} w + s)\frac{\partial }{\partial z_2 }+}
\hfill \\
\qquad
 {+(c_{1} z_1 + c_{2} z_2 +c_{3} w  + q)\frac{\partial }{\partial w
 }}
\hfill \\
\end{array} }
\eqno (2.11)
$$
a matrix of the form
$$
    Z =
 \left(
\begin{array}{cccc}
a_{1} & a_{2} & a_{3} & p\\
b_{1} & b_{2} & b_{3} & s\\
c_{1} & c_{2} & c_{3} & q\\
  0    &   0    &   0    & 0
\end{array}
\right). \eqno (2.12)
$$

   Here the Lie bracket of two vector fields
corresponds to the usual commutator for a matrices
$$
 [Z_1, Z_2] = Z_1 Z_2 - Z_2 Z_1.
$$

  Let $ g(M) $ be a matrix Lie algebra related to the
affine homogeneous hypersurface $ M \subset \mathbb{C}^3 $ and the
surface itself is defined by canonical equation (1.1). The fourth
column of the matrix (2.12) corresponds to the translation component
of vector field.  That's why the real linear span of the fourth
columns of the matrices (2.12) coincides with the tangent plane $
T_0 M $ and has dimension 5. Hence, $ \dim_{\mathbb{R}} g(M) \ge 5
$.

  Starting from the part 3 of this article, we consider only 5-dimensional Lie algebras with
  the bases ($ m_1 $, $ m_2 $, $ m_3 $, $ m_4 \in \mathbb{R} $)
$$
E_{1}= \left(
  \begin{array}{cccc}
    A1_1 & A2_1 & A3_1 & 1 \\
    B1_1 & B2_1 & B3_1 & 0 \\
    4i & 0 & m_1 & 0 \\
    0 & 0 & 0 & 0 \\
  \end{array}
\right), E_{2}= \left(
  \begin{array}{cccc}
    A1_2 & A2_2 & A3_2 & i \\
    B1_2 & B2_2 & B3_2 & 0 \\
    0 & 0 & m_2 & 0 \\
    0 & 0 & 0 & 0 \\
  \end{array}
\right),$$ $$E_{3}= \left(
  \begin{array}{cccc}
    A1_3 & A2_3 & A3_3 & 0 \\
    B1_3 & B2_3 & B3_3 & 1 \\
    0 & 2i & m_3 & 0 \\
    0 & 0 & 0 & 0 \\
  \end{array}
\right),\eqno (2.13)
$$
$$
E_{4}= \left(
  \begin{array}{cccc}
    A1_4 & A2_4 & A3_4 & 0 \\
    B1_4 & B2_4 & B3_4 & i \\
    0 & 2 & m_4 & 0 \\
    0 & 0 & 0 & 0 \\
  \end{array}
\right), E_{5}= \left(
  \begin{array}{cccc}
    A1_5 & A2_5 & A3_5 & 0 \\
    B1_5 & B2_5 & B3_5 & 0 \\
    0 & 0 & m_5 & 1 \\
    0 & 0 & 0 & 0 \\
  \end{array}
\right).
$$

  It can be shown (see [13], [14]) that for the elements $c_1, c_2$ of
the matrices (2.12) the following restrictions hold
$$
  q \in \mathbb{R}, \quad
    c_{1} = 2i (\bar p + 2\varepsilon_1 p), \quad
    c_{2} = 2i (\bar s + 2\varepsilon_1 s).
$$

  We declare that for the hypersurfaces of the type (1/2, \,0) it is
sufficient to consider only 5-dimensional algebras and that all
algebras of higher dimensions can be reduced to 5-dimensional ones.
   Here we omit the proof of this statement but it will be a part of a
future big article that authors prepare.  We only add that for
algebras $g(M)$ corresponding to homogeneous (1/2, \,0)-type
surfaces (similarly to the surfaces from the paper [14]), the
following upper dimension estimate holds
$$
  \dim_{\mathbb{R}} g(M) \le 7.
\eqno (2.14)
$$

   The equality
$
  \dim_{\mathbb{R}} g(M) = 7
$ is satisfied only for unique (up to affine equivalence) affine
homogeneous hypersurface of the type (1/2, \,0), namely for the
quadric
$$
   v = |z_1|^2 + |z_2|^2 + \frac12 (z_1^2 + \bar z_1^2) = 2 x_1^2 + |z_2|^2.
\eqno (2.15)
$$

  As a basis of the algebra of affine vector fields for this surface we can use,
for example, the matrices
$$
E_{1}= \left(
  \begin{array}{cccc}
    0 & 0 & 0 & 1 \\
    0 & 0 & 0 & 0 \\
    4i & 0 & 0 & 0 \\
    0 & 0 & 0 & 0 \\
  \end{array}
\right), E_{2}= \left(
  \begin{array}{cccc}
    0 & 0 & 0 & i \\
    0 & 0 & 0 & 0 \\
    0 & 0 & 0 & 0 \\
    0 & 0 & 0 & 0 \\
  \end{array}
\right), E_{3}= \left(
  \begin{array}{cccc}
    0 & 0 & 0 & 0 \\
    0 & 0 & 0 & 1 \\
    0 & 2i & 0 & 0 \\
    0 & 0 & 0 & 0 \\
  \end{array}
\right),
$$
$$
E_{4}= \left(
  \begin{array}{cccc}
    0 & 0 & 0 & 0 \\
    0 & 0 & 0 & i \\
    0 & 2 & 0 & 0 \\
    0 & 0 & 0 & 0 \\
  \end{array}
\right), E_{5}= \left(
  \begin{array}{cccc}
    0 & 0 & 0 & 0 \\
    0 & 0 & 0 & 0 \\
    0 & 0 & 0 & 1 \\
    0 & 0 & 0 & 0 \\
  \end{array}
\right), \eqno (2.16)
$$
$$
E_{6}= \left(
    \begin{array}{cccc}
      1 & 0 & 0 & 0 \\
      0 & 1 & 0 & 0 \\
      0 & 0 & 2 & 0 \\
      0 & 0 & 0 & 0 \\
    \end{array}
  \right),
E_{7}= \left(
    \begin{array}{cccc}
      0 & 0 & 0 & 0 \\
      0 & i & 0 & 0 \\
      0 & 0 & 0 & 0 \\
      0 & 0 & 0 & 0 \\
    \end{array}
  \right).
$$

\

   The algebra with the basis (2.16) will be used below for simplification of
cumbersome integration of some 5-dimensional algebras. However, the
main difficulty in our article is to construct a complete list of
such algebras related to affine homogeneous hypersurfaces.

To obtain this list one need to study a system of quadratic
equations which means the bracket-closedness of the linear space
with the basis (2.13). Element-wise rewriting of 10 such matrix
equations gives us the system of 120 scalar ones.

Using the specific form of the basis (2.13) one can extract (like in
[14]) some relatively simple subsystems from this big system. Their
solutions can be obtained by consideration of a large number of
partial cases using any computer algebra system. On the final step
of the described procedure we obtain the desired list of algebras.
  One can mention here [3] as one of the first papers studying the Lie
algebras and homogeneity by means of computer mathematics. Note also
the recent work [19] where the main results were obtained in the
same way.

\

{\bf Remark.} Note that complete classification of all the
5-dimensional algebras (see [18]) is well-known. However, we don't
know the simple method to choose algebras related to homogeneous
hypersurfaces from this complete list of algebras.

\

\centerline{\textbf{3. Lie algebras related to homogeneous
hypersurfaces of a type (1/2, \,0)}}

\

  Here we give the bases of nine types of 5-dimensional algebras
related to homogeneous hypersurfaces, that were obtained by the
method mentioned above. All the parameters $ m_k, t_j $ using below
are arbitrary real numbers.

\
$$
E_1 = \left(%
\begin{array}{cccc}
  m_1 & 0 & 0 & 1 \\
  0 & m_1+\frac{im_1t_{16}}{m_2} & 0 & 0 \\
  4i & 0 & 2m_1 & 0 \\
  0 & 0 & 0 & 0 \\
\end{array}%
\right),\eqno (3.1)$$ $$E_2 = \left(%
\begin{array}{cccc}
  m_2 & 0 & 0 & i \\
  0 & m_2+it_{16} & 0 & 0 \\
  0 & 0 & 2m_2 & 0 \\
  0 & 0 & 0 & 0 \\
\end{array}%
\right),
$$
$$
E_3 = \left(%
\begin{array}{cccc}
  0 & 0 & 0 & 0 \\
  0 & 0 & 0 & 1 \\
  0 & 2i & 0 & 0 \\
  0 & 0 & 0 & 0 \\
\end{array}%
\right), E_4 = \left(%
\begin{array}{cccc}
  0 & 0 & 0 & 0 \\
  0 & 0 & 0 & i \\
  0 & 2 & 0 & 0 \\
  0 & 0 & 0 & 0 \\
\end{array}%
\right), \
E_5 = \left(%
\begin{array}{cccc}
  0 & 0 & 0 & 0 \\
  0 & 0 & 0 & 0 \\
  0 & 0 & 0 & 1 \\
  0 & 0 & 0 & 0 \\
\end{array}%
\right).
$$

\

\

$$
 E_1 = \left(%
\begin{array}{cccc}
  2t_7 - 2m_1 - 2im_2 & 0 & i(m_1+im_2)(m_1-t_7) & 1 \\
  0 & t_7 + it_8 & 0 & 0 \\
  4i & 0 & 2m_1 & 0 \\
  0 & 0 & 0 & 0 \\
\end{array}%
\right), \eqno (3.2)
$$
$$
E_2 = \left(%
\begin{array}{cccc}
  2m_2 & 0 & 0 & i \\
  0 & m_2+it_{16} & 0 & 0 \\
  0 & 0 & 2m_2 & 0 \\
  0 & 0 & 0 & 0 \\
\end{array}%
\right),
$$
$$
E_3 = \left(%
\begin{array}{cccc}
  0 & -m_1 + t_7 & 0 & 0 \\
  0 & 0 & 0 & 1 \\
  0 & 2i & 0 & 0 \\
  0 & 0 & 0 & 0 \\
\end{array}%
\right), E_4 = \left(%
\begin{array}{cccc}
  0 & i(m_1 - t_7) & 0 & 0 \\
  0 & 0 & 0 & i \\
  0 & 2 & 0 & 0 \\
  0 & 0 & 0 & 0 \\
\end{array}%
\right),
$$
$$
E_5 = \left(%
\begin{array}{cccc}
  -m_2(m_1-t_7) & 0 & 0 & 0 \\
  0 & -\frac{1}{2}(m_1 - t_7)(m_2+it_{16}) & 0 & 0 \\
  0 & 0 & -m_2(m_1-t_7) & 1 \\
  0 & 0 & 0 & 0 \\
\end{array}%
\right).$$

\

\

$$
E_1 = \left(%
\begin{array}{cccc}
  t_1 & 0 & -\frac{i}{2}(m_1-t_7)(t_1 - 2 t_7) & 1 \\
  0 & t_7 +i t_8 & 0 & 0 \\
  4i & 0 & 2m_1 & 0 \\
  0 & 0 & 0 & 0 \\
\end{array}%
\right),\eqno (3.3)$$

$$E_2 = \left(%
\begin{array}{cccc}
  0 & 0 & 0 & i \\
  0 & 0 & 0 & 0 \\
  0 & 0 & 0 & 0 \\
  0 & 0 & 0 & 0 \\
\end{array}%
\right), \;\;
E_3 = \left(%
\begin{array}{cccc}
  0 & -(m_1 - t_7) & 0 & 0 \\
  0 & 0 & 0 & 1 \\
  0 & 2i & 0 & 0 \\
  0 & 0 & 0 & 0 \\
\end{array}%
\right), $$
$$E_4 = \left(%
\begin{array}{cccc}
  0 & i(m_1 - t_7) & 0 & 0 \\
  0 & 0 & 0 & i \\
  0 & 2 & 0 & 0 \\
  0 & 0 & 0 & 0 \\
\end{array}%
\right),
E_5 = \left(%
\begin{array}{cccc}
  0 & 0 & 0 & 0 \\
  0 & 0 & 0 & 0 \\
  0 & 0 & 0 & 1 \\
  0 & 0 & 0 & 0 \\
\end{array}%
\right).
$$

\

\

$$
E_1 = \left(%
\begin{array}{cccc}
  -2m_1 - 2im_2 & 0 & im_1(m_1+im_2) & 1 \\
  0 & 2im_2 & 0 & 0 \\
  4i & 0 & 2m_1 & 0 \\
  0 & 0 & 0 & 0 \\
\end{array}%
\right), \eqno (3.4)$$ $$E_2 = \left(%
\begin{array}{cccc}
  2m_2 & 0 & 0 & i \\
  0 & 2m_2 + it_{16} & 0 & 0 \\
  0 & 0 & 2m_2 & 0 \\
  0 & 0 & 0 & 0 \\
\end{array}%
\right),
E_3 = \left(%
\begin{array}{cccc}
  0 & -m_1-im_2 & 0 & 0 \\
  -2im_2 & 0 & -m_1m_2 & 1 \\
  0 & 2i & 0 & 0 \\
  0 & 0 & 0 & 0 \\
\end{array}%
\right),$$ $$E_4 = \left(%
\begin{array}{cccc}
  0 & im_1-m_2 & 0 & 0 \\
  2m_2 & 0 & -im_1m_2 & i \\
  0 & 2 & 0 & 0 \\
  0 & 0 & 0 & 0 \\
\end{array}%
\right),
$$
$$
E_5 = \left(%
\begin{array}{cccc}
  -m_1m_2 & 0 & 0 & 0 \\
  0 & -\frac{1}{2}m_1(2m_2+it_{16}) & 0 & 0 \\
  0 & 0 & -m_1m_2 & 1 \\
  0 & 0 & 0 & 0 \\
\end{array}%
\right).$$

\

\

$$
E_1 = \left(%
\begin{array}{cccc}
  t_1 & 2m_3-2im_4 & -\frac{i}{2}(t_1m_1+2m_3^2+2m_4^2) & 1 \\
  0 & 0 & 0 & 0 \\
  4i & 0 & 2m1 & 0 \\
  0 & 0 & 0 & 0 \\
\end{array}%
\right),\eqno (3.5)$$ $$E_2 = \left(%
\begin{array}{cccc}
  0 & 0 & 0 & i \\
  0 & 0 & 0 & 0 \\
  0 & 0 & 0 & 0 \\
  0 & 0 & 0 & 0 \\
\end{array}%
\right),
E_3 = \left(%
\begin{array}{cccc}
  2m_3 & -m_1 & 0 & 0 \\
  0 & m_3-im_4 & 0 & 1 \\
  0 & 2i & 2m_3 & 0 \\
  0 & 0 & 0 & 0 \\
\end{array}%
\right),$$ $$E_4 = \left(%
\begin{array}{cccc}
  2m_4 & im_1 & 0 & 0 \\
  0 & m_4+im_3 & 0 & i \\
  0 & 2 & 2m_4 & 0 \\
  0 & 0 & 0 & 0 \\
\end{array}%
\right), E_5 = \left(%
\begin{array}{cccc}
  0 & 0 & 0 & 0 \\
  0 & 0 & 0 & 0 \\
  0 & 0 & 0 & 1 \\
  0 & 0 & 0 & 0 \\
\end{array}%
\right).
$$

\

\

$$
E_1 = \left(%
\begin{array}{cccc}
  -2m_1 & m_3-im_4 & im_1^2 & 1 \\
  0 & 0 & 0 & 0 \\
  4i & 0 & 2m_1 & 0 \\
  0 & 0 & 0 & 0 \\
\end{array}%
\right), E_2 = \left(%
\begin{array}{cccc}
  0 & 0 & 0 & i \\
  0 & 0 & 0 & 0 \\
  0 & 0 & 0 & 0 \\
  0 & 0 & 0 & 0 \\
\end{array}%
\right), \eqno (3.6)
$$
$$
E_3 = \left(%
\begin{array}{cccc}
  \frac{3}{2}m_3 & -m_1 & \frac{i}{4}m_1m_3 & 0 \\
  0 & m_3 & 0 & 1 \\
  0 & 2i & 2m_3 & 0 \\
  0 & 0 & 0 & 0 \\
\end{array}%
\right),
$$
$$
E_4 = \left(%
\begin{array}{cccc}
  \frac{3}{2}m_4 & im_1 & \frac{i}{4}m_1m_4 & 0 \\
  0 & m_4 & 0 & i \\
  0 & 2 & 2m_4 & 0 \\
  0 & 0 & 0 & 0 \\
\end{array}%
\right), E_5 = \left(%
\begin{array}{cccc}
  0 & 0 & 0 & 0 \\
  0 & 0 & 0 & 0 \\
  0 & 0 & 0 & 1 \\
  0 & 0 & 0 & 0 \\
\end{array}%
\right).
$$

\

\

$$
E_1 = \left(%
\begin{array}{cccc}
  t_1 & t_3+it_4 & -\frac{i}{2}(t_1m_1+t_4^2+t_3^2) & 1 \\
  2t_3-2it_4 & 0 & im_1(-t_3+it_4) & 0 \\
  4i & 0 & 2m_1 & 0 \\
  0 & 0 & 0 & 0 \\
\end{array}%
\right), \eqno (3.7)$$
$$E_2 = \left(%
\begin{array}{cccc}
  0 & t_4-it_3 & \frac{1}{2}(t_3^2+t_4^2) & i \\
  0 & 0 & 0 & 0 \\
  0 & 0 & 0 & 0 \\
  0 & 0 & 0 & 0 \\
\end{array}%
\right),
$$
$$
E_3 = \left(%
\begin{array}{cccc}
  t_3-it_4 & -m_1 & \frac{i}{2}m_1(-t_3+it_4) & 0 \\
  0 & -2it_4 & \frac{i}{2}(-t_3+it_4)(t_3+it_4) & 1 \\
  0 & 2i & 0 & 0 \\
  0 & 0 & 0 & 0 \\
\end{array}%
\right),
$$
$$
E_4 = \left(%
\begin{array}{cccc}
  -t_4-it_3 & im_1 & \frac{1}{2}m_1(-t_3+it_4) & 0 \\
  0 & -2it_3 & \frac{1}{2}(t_3^2 + t_4^2) & i \\
  0 & 2 & 0 & 0 \\
  0 & 0 & 0 & 0 \\
\end{array}%
\right),
$$
$$
E_5 = \left(%
\begin{array}{cccc}
  -\frac{i}{2}(t_3^2+t_4^2) & \frac{i}{2}m_1(t_3+it_4) & -\frac{1}{4}m_1(t_3^2+t_4^2) & 0 \\
  0 & -\frac{i}{2}(t_3^2+t_4^2) & 0 & 0 \\
  0 & 0 & -\frac{i}{2}(t_3^2+t_4^2) & 1 \\
  0 & 0 & 0 & 0 \\
\end{array}%
\right).
$$

\

\

$$
E_1 = \left(%
\begin{array}{cccc}
  -2m_1 & 0 & im_1^2 & 1 \\
  0 & 0 & 0 & 0 \\
  4i & 0 & 2m_1 & 0 \\
  0 & 0 & 0 & 0 \\
\end{array}%
\right), E_2 = \left(%
\begin{array}{cccc}
  0 & 0 & 0 & i \\
  0 & 0 & 0 & 0 \\
  0 & 0 & 0 & 0 \\
  0 & 0 & 0 & 0 \\
\end{array}%
\right), \eqno (3.8)
$$
$$
E_3 = \left(%
\begin{array}{cccc}
  m_3 & -m_1 & \frac{i}{2}m_1m_3 & 0 \\
  0 & m_3-im_4 & 0 & 1 \\
  0 & 2i & 2m_3 & 0 \\
  0 & 0 & 0 & 0 \\
\end{array}%
\right),
$$
$$
E_4 = \left(%
\begin{array}{cccc}
  m_4 & im_1 & \frac{i}{2}m_1m_4 & 0 \\
  0 & m_4+im_3 & 0 & i \\
  0 & 2 & 2m_4 & 0 \\
  0 & 0 & 0 & 0 \\
\end{array}%
\right), E_5 = \left(%
\begin{array}{cccc}
  0 & 0 & 0 & 0 \\
  0 & 0 & 0 & 0 \\
  0 & 0 & 0 & 1 \\
  0 & 0 & 0 & 0 \\
\end{array}%
\right).
$$

\

\

$$
E_1 = \left(%
\begin{array}{cccc}
  3t_7-2m_1 & 0 & \frac{i}{2}(2m_1-t_7)(m_1-t_7) & 1 \\
  0 & t_7 & 0 & 0 \\
  4i & 0 & 2m_1 & 0 \\
  0 & 0 & 0 & 0 \\
\end{array}%
\right), \eqno (3.9)$$ $$E_2 = \left(%
\begin{array}{cccc}
  m_2 & 0 & \frac{i}{2}m_2(m_1-t_7) & i \\
  0 & m_2 & 0 & 0 \\
  0 & 0 & 2m_2 & 0 \\
  0 & 0 & 0 & 0 \\
\end{array}%
\right),
E_3 = \left(%
\begin{array}{cccc}
  m_3 & t_7-m_1 & \frac{i}{2}m_3(m_1-t_7) & 0 \\
  0 & m_3 & 0 & 1 \\
  0 & 2i & 2m_3 & 0 \\
  0 & 0 & 0 & 0 \\
\end{array}%
\right),$$ $$E_4 = \left(%
\begin{array}{cccc}
  m_4 & i(m_1-t_7) & \frac{i}{2}m_4(m_1-t_7) & 0 \\
  0 & m_4 & 0 & i \\
  0 & 2 & 2m_4 & 0 \\
  0 & 0 & 0 & 0 \\
\end{array}%
\right),
$$
$$
E_5 = \left(%
\begin{array}{cccc}
  -\frac{1}{2}m_2(m_1-t_7) & 0 & -\frac{i}{4}m_2(m_1-t_7)^2 & 0 \\
  0 & -\frac{1}{2}m_2(m_1-t_7) & 0 & 0 \\
  0 & 0 & -m_2(m_1-t_7) & 1 \\
  0 & 0 & 0 & 0 \\
\end{array}%
\right).
$$

\

\centerline{\textbf{4. Integration of Lie algebras}}

\

   Note that integrating of a particular Lie algebra can give
different (affine inequivalent) homogeneous hypersurfaces. Hence,
the classifications results for homogeneous manifolds that were
formulated in terms of Lie algebras (see, for example [2]) need to
be specified.

  In this section we give the results of integration of
nine types of algebras in expanded form. As it was mentioned above,
integration of algebras (i.e. constructing the homogeneous
hypersurfaces related to these algebras) is a complicated procedure.
So in the section 5 we outline shortly only some details of
integration of the algebras under study.

\

{\bf Theorem 4.1.} {\it Integral variety of a type $(1/2,\,0)$
related to any of the algebras (3.1), (3.8), (3.9) is affine
equivalent to a quadric
$$
  v = 2 x_1^2 + |z_2|^2.
$$
}

{\bf Theorem 4.2.} {\it Integral variety of a type $(1/2,\,0)$
related to any of the algebras (3.2) is affine equivalent to one of
the following hypersurfaces:
$$
  1) \quad
   v = \exp (x_1) + |z_2|^2,
$$
$$
  2) \quad   v = |z_1|e^{B \arg{z_1}} + |z_2|^2, \ B \in \mathbb{R};
$$
}

{\bf Theorem 4.3.} {\it Integral variety of a type $(1/2,\,0)$
related to any of the algebras (3.3) is affine equivalent to one of
the following hypersurfaces:
$$
  1) \quad
      v = 2 x_1^2 + |z_2|^2,
$$
$$
  2) \quad
   v = \exp (x_1) + |z_2|^2,
$$
$$
  3) \quad
   v = - \ln(1 + x_1) + |z_2|^2,
$$
$$
  4) \quad
   v = (1+x_1) \ln(1 + x_1) + |z_2|^2,
$$
$$
  5) \quad   v = \pm (1+x_1)^{\alpha} + |z_2|^2, \quad
       \alpha \in \mathbb{R} \setminus\{0,1,2\}.
$$
}

{\bf Theorem 4.4.} {{\it Integral variety of a type $(1/2,\,0)$
related to any of the algebras (3.4) is affine equivalent to one of
the following hypersurfaces:
$$
  1) \quad
      v = 2 x_1^2 + |z_2|^2,
$$
$$
  2) \quad
   v^2 = |z_1|^2 + |z_2|^2.
$$
}

{\bf Theorem 4.5.} {\it Integral variety of a type $(1/2,\,0)$
related to any of the algebras (3.5) is affine equivalent to one of
the following hypersurfaces:
$$
  1) \quad
      v = 2 x_1^2 + |z_2|^2,
$$
$$
  2) \quad
   v = (1+x_1) \ln(1 + x_1) + |z_2|^2,
$$
$$
  3) \quad
   v = x_1^{(1-\alpha)}|z_2|^{2\alpha}, \ \alpha \in \mathbb{R} \setminus \{0,1\}.
$$
}

{\bf Theorem 4.6.} {\it Integral variety of a type $(1/2,\,0)$
related to any of the algebras (3.6) is affine equivalent to one of
the following hypersurfaces:
$$
   1) \ v = 2 x_1^2 + |z_2|^2,
$$
$$
   2) \ v = \frac{x_1^2}{1 - x_2} + 2 |z_2|^2.
$$
}

{\bf Theorem 4.7} {\it Integral variety of a type $(1/2,\,0)$
related to any of the algebras (3.7) is affine equivalent to one of
the following hypersurfaces:
$$
   1) \quad  v = 2x_1^2 + |z_2|^2,
$$
$$
   2) \quad  v = -\ln(1+ x_1) + |z_2|^2,
$$
$$
   3) \
Re\,(\bar z_1 w) = (Re\,(z_1 \bar z_2))^{\alpha}, \quad \alpha \in
(- \infty,0).
$$
}
  It is obvious that combination of these seven theorems leads to
the main theorem stated above.

\

\centerline{\textbf{5. Remarks on the proofs of the theorems }}

\

  First of all we note the role played by 7-dimensional Lie algebra (2.16)
in the following reasonings. This algebra corresponds to the quadric
$
    v = 2 x_1^2 + |z_2|^2.
$
  Many of the algebras mentioned above can be reduced to 5-dimensional
subalgebras of (2.16) by means of suitable matrix similarity. This
argument leads to the simple proof of the theorem 4.1.

\

\centerline{\textbf{5.1 Proof of the theorem 4.1}}

\

  Any algebra in the family (3.1) is just a subalgebra of the algebra (2.16).
This means that integral manifolds of such subalgebras coincide with
the analogous manifold related to (2.16), i.e. any of them is a
quadric (2.15).

  Algebras in the families (3.8) and (3.9) are reduced to the subalgebras
of (2.16) by means of matrix similarities. For instance, in the case
(3.8) this is the similarity $
    g^* = C^{-1} g C,
$ with the matrix
$$
C =
\left(%
\begin{array}{cccc}
  1 & 0 & i\lambda & 0 \\
  0 & 1 & 0 & 0 \\
  0 & 0 & 1 & 0 \\
  0 & 0 & 0 & 1 \\
\end{array}%
\right), \eqno (5.1)
$$
where $ \lambda = m_1/2. $

  A new basis of (3.8) takes the form
$$
E_{1}= \left(
  \begin{array}{cccc}
    0 & 0 & 0 & 1 \\
    0 & 0 & 0 & 0 \\
    4i & 0 & 0 & 0 \\
    0 & 0 & 0 & 0 \\
  \end{array}
\right), E_{2}= \left(
  \begin{array}{cccc}
    0 & 0 & 0 & i \\
    0 & 0 & 0 & 0 \\
    0 & 0 & 0 & 0 \\
    0 & 0 & 0 & 0 \\
  \end{array}
\right), \eqno (5.2)$$ $$E_{3}= \left(
  \begin{array}{cccc}
    m_3 & 0 & 0 & 0 \\
    0 & m_3-im_4 & 0 & 1 \\
    0 & 2i & 2m_3 & 0 \\
    0 & 0 & 0 & 0 \\
  \end{array}
\right),
$$
$$
E_{4}= \left(
  \begin{array}{cccc}
    m_4 & 0 & 0 & 0 \\
    0 & m_4+im_3 & 0 & i \\
    0 & 2 & 2m_4 & 0 \\
    0 & 0 & 0 & 0 \\
  \end{array}
\right), E_{5}= \left(
  \begin{array}{cccc}
    0 & 0 & 0 & 0 \\
    0 & 0 & 0 & 0 \\
    0 & 0 & 0 & 1 \\
    0 & 0 & 0 & 0 \\
  \end{array}
\right).
$$

\

  Analogously any algebra (3.9) also turns
under the similarity with matrix (5.1) and $ \lambda = (m_1 - t_7)/2
$ into subalgebra of (2.16).

  Theorem 4.1 is proved.

\

  The scheme of the study of another algebras is more complicated. Here
one need to work with the systems of partial differential equations.

  The fact that the vector field $ Z $ is tangent to the homogeneous surface
$
   M = \{ \Phi(z, \bar z, u, v) = 0   \}
$ can be expressed in the form
$$
   Re \left(Z(\Phi)|_{M}\right) = 0.
\eqno (5.3)
$$

  So we have the system of five equations (5.3)
related to the basis vector fields for each algebra under
consideration.
  The defining function
$$
  \Phi(z, \bar z, u, v) = - v + F(z, \bar z, u)
$$
of homogeneous surface determined by unknown function $ F(z, \bar z,
u) $ can be found from every such system (5.3).

\

\centerline{\textbf{5.2 Discussion of the proof of the theorem 4.2}}

\

  As in the previous case, here we use the similarity with the matrix (5.1)
and the same $ \lambda = (m_1 - t_7)/2 $. But we consider the sum $
   E_5 + \frac{m_1 - t_7}2 E_2
$ instead of the matrix $ E_5 $.

    Then the basis of the new algebra takes the form
$$
E_{1}= \left(
  \begin{array}{cccc}
    2im_2 & 0 & 0 & 1 \\
    0 & t_7+it_8 & 0 & 0 \\
    4i & 0 & 2t_7 & 0 \\
    0 & 0 & 0 & 0 \\
  \end{array}
\right),\eqno (5.4)$$ $$E_{2}= \left(
  \begin{array}{cccc}
    m_2 & 0 & 0 & i \\
    0 & m_2 + it_{16} & 0 & 0 \\
    0 & 0 & 2m_2 & 0 \\
    0 & 0 & 0 & 0 \\
  \end{array}
\right), E_{3}= \left(
  \begin{array}{cccc}
    0 & 0 & 0 & 0 \\
    0 & 0 & 0 & 1 \\
    0 & 2i & 0 & 0 \\
    0 & 0 & 0 & 0 \\
  \end{array}
\right),
$$ $$E_{4}= \left(
  \begin{array}{cccc}
    0 & 0 & 0 & 0 \\
    0 & 0 & 0 & i \\
    0 & 2 & 0 & 0 \\
    0 & 0 & 0 & 0 \\
  \end{array}
\right), E_{5}= \left(
  \begin{array}{cccc}
    0 & 0 & 0 & 0 \\
    0 & 0 & 0 & 0 \\
    0 & 0 & 0 & 1 \\
    0 & 0 & 0 & 0 \\
  \end{array}
\right).
$$

   Note that one can easily integrate this algebra
due to the presence of the triple of "trivial" \ matrices $ E_3,
E_4, E_5 $ in its basis. For instance, the equation (5.3) related to
the field $ E_5 $, means the independence of the function $ F(z,
\bar z, u) $ on the variable $ u = Re\, w $ (the {\it rigidity} of
the surface under consideration in terms of several dimensional
complex analysis) in this case.

  Another two simple equations are related to the fields $ E_3, E_4 $
$$
  Re \left(\frac{\partial F}{\partial z_2} + 2\, i\, z_2\,\frac{i}2 \right) = 0,
\qquad
  Re \left(i\frac{\partial F}{\partial z_2} + 2\, z_2\, \frac{i}2 \right) = 0.
$$

   Hence,
$$
  \frac{\partial F}{\partial z_2}  = \bar z_2, \quad
\mbox{i.e.} \quad F = F(z,\bar z) = |z_2|^2 + H(z_1, \bar z_1),
\eqno (5.5)
$$
where $ H(z_1, \bar z_1) $ is an arbitrary real-valued function on
the variables $ z_1,\bar z_1 $.

  Two remaining more complicated equations (5.3) related to the fields
$ E_1, E_2 $ of the basis (5.2) have the form

$$
    (2\,m_2 y_1 + 1) \frac{\partial H}{\partial x_1}
       - 2\,m_2 x_1 \frac{\partial H}{\partial y_1} =
          2 t_7 H,
\eqno (5.6)
$$
$$
    2\,m_2 x_1 \frac{\partial H}{\partial x_1}
     + (2\,m_2 y_1 + 1) \frac{\partial H}{\partial y_1} = 2 m_2 \,H.
$$

   Solving of this system of equations leads to two cases
associated with the possible values of the parameter $ m_2 $.

   In the first case when $ m_2 = 0 $, the general solution of the system
(5.6) has the form
$$
   H = C \exp (2 t_7 x_1)
\eqno (5.7)
$$
with arbitrary constant $ C $.

   It remains to note that in the case $ C = 0 $ (as well as in the case
$ C \ne 0, t_7 = 0 $) formula (5.7) defines a Levi degenerate
surface, while the condition $ C < 0 $ leads to the surface with a
non-degenerate indefinite Levi form. We study only strictly
pseudo-convex homogeneous surfaces, so it is necessary to consider
only positive values of the constant $ C $ and non-zero values of
the parameter $ t_7 $ in formula (5.7). In this case the equation of
the required affine-homogeneous surface corresponding to (5.7) can
be written as
$$
   v = \exp (2 t_7 x_1 + \ln C) + |z_2|^2,
$$
or in the form (2.3) after the change of variable $
  x_1^* =  2 t_7 x_1 + \ln C.
$

   In the second case, i.e. if $ m_2 \ne 0 $, the formal solution of
the system (5.6) depends on two real parameters. However, by
reasonings similar to those of the case $m_2 = 0$, we can rewrite
the solution in the form (2.8) so it depends on just one parameter.

\

\centerline{\textbf{5.3 Comments to the proof of the Theorem 4.3}}

\

   Using the similarity with the matrix (5.1) and taking
$ \lambda = (m_1 - t_7)/2 $ we can construct new algebras with bases
$$
  E_1 =  \left( \begin {array}{cccc} r &0&0&1
\\\noalign{\medskip}0&t_{{7}}+it_{{8}}&0&0\\\noalign{\medskip}4\,i&0&2
\,t_{{7}}&0\\\noalign{\medskip}0&0&0&0\end {array} \right), E_{2}=
\left(
  \begin{array}{cccc}
    0 & 0 & 0 & i \\
    0 & 0 & 0 & 0 \\
    0 & 0 & 0 & 0 \\
    0 & 0 & 0 & 0 \\
  \end{array}
\right),\eqno (5.8)$$ $$E_{3}= \left(
  \begin{array}{cccc}
    0 & 0 & 0 & 0 \\
    0 & 0 & 0 & 1 \\
    0 & 2i & 0 & 0 \\
    0 & 0 & 0 & 0 \\
  \end{array}
\right), E_{4}= \left(
  \begin{array}{cccc}
    0 & 0 & 0 & 0 \\
    0 & 0 & 0 & i \\
    0 & 2 & 0 & 0 \\
    0 & 0 & 0 & 0 \\
  \end{array}
\right), E_{5}= \left(
  \begin{array}{cccc}
    0 & 0 & 0 & 0 \\
    0 & 0 & 0 & 0 \\
    0 & 0 & 0 & 1 \\
    0 & 0 & 0 & 0 \\
  \end{array}
\right),
$$
where $r = t_{{1}} + 2\,m_{{1}} - 2\,t_{{7}}$.

   By analogy with the arguments above, we concentrate our attention on the field $ E_1 $. After
straightforward simplifications associated with the other basic
fields, the equation of a surface under consideration takes the form
$$
    F = |z_2|^2 + H(x_1).
\eqno (5.9)
$$

  Here $ H(x_1) $  is an analytic function presenting a solution of ODE
$$
   (1+ r x_1) H'(x_1) = 2t_7 H + 4 x_1.
\eqno (5.10)
$$

   The form of solution of (5.10) depends on the parameters $t_7$ and $r$.
 Thus, when
$ t_7 = 0, r = 0 $ we get
$$
   H(x_1) = 2 x_1^2 + C,
$$
and the equation of homogeneous surface has the form (2.1). If $ t_7
= 0, r \ne 0 $, the equation (2.2) appears. If $ t_7 \ne 0 $, but $
r = 0 $, we get the equation (2.3).

   Finally, in the general case, when $ t_7 \ne 0 $, $ r \ne 0 $, there are two
types of solutions of (5.10):
   when $ {2t_7} / r = 1 $ we have the equation (2.5) of the
homogeneous surface and if $ \alpha = {2t_7} / r \ne 1 $ we obtain
the equation (2.4).

\

\centerline{\textbf{5.4 The main points of the proof of the Theorem
4.4}}

\

  Here, as in the previous case, it is convenient to start with the similarity
transformation with the matrix (5.1) where $ \lambda = m_1 / 2 $. As
a result, the parameter $ m_1 $ is excluded from the basis of the
algebra.

    Thereby, the parameter $ m_2 $ will play the main role in the
discussion. Value $ m_2 = 0 $ yields a family of 5-dimensional
algebras, each of which is a subalgebra of (2.16). So the quadric $
  v = 2 x_1^2 + |z_2|^2
$ appears in Theorem 4.4.

  If $ m_2 \neq 0 $, any algebra of the family (3.4) is similar to the
algebra with basis
$$
E_1 =
  \left( \begin {array}{cccc} -i&0&0&0\\\noalign{\medskip}0&i&0&0
\\\noalign{\medskip}0&0&0&0\\\noalign{\medskip}0&0&0&0\end {array}
 \right)
, \ E_2 =
 \left( \begin {array}{cccc} 1&0&0&0\\\noalign{\medskip}0&1+i\xi&0&0
\\\noalign{\medskip}0&0&1&0\\\noalign{\medskip}0&0&0&0\end {array}
 \right)
,\eqno (5.11)$$ $$\ E_3 =
 \left( \begin {array}{cccc} 0&i&0&0\\\noalign{\medskip}i&0&0&0
\\\noalign{\medskip}0&0&0&0\\\noalign{\medskip}0&0&0&0\end {array}
 \right)
, E_4 =
 \left( \begin {array}{cccc} 0&1&0&0\\\noalign{\medskip}-1&0&0&0
\\\noalign{\medskip}0&0&0&0\\\noalign{\medskip}0&0&0&0\end {array}
 \right)
, \ E_5 =
 \left( \begin {array}{cccc} 0&0&0&0\\\noalign{\medskip}0&0&0&0
\\\noalign{\medskip}0&0&0&1\\\noalign{\medskip}0&0&0&0
\end {array}
\right) ,
$$
where $ \xi = t_{16}/2m_2 \in \mathbb{R} $.

  One remark concerning the upper left $ 2 \times 2 $ matrix block of
the matrices $ E_1, E_3, E_4 $ (responsible for the action of the
group $ G(M) $ in the complex space $ \mathbb{C}^2_{(z_1, z_2)}$
tangent to discussed homogeneous surfaces) essentially helps in
building an integral manifold of this algebra. These three $ 2
\times 2 $ Pauli matrices present the basis of the well-known
algebra $ su(2) $. It is a Lie algebra for the group $ SU (2) $
which indicates the occurrence of the expression $ | z_1 | ^ 2 + |
z_2 | ^ 2 $ in equation of discussed homogeneous hypersurfaces.
   Thus the appearance of the equation
(2.6) of affine-homogeneous surface can be naturally explained in
this case.

   Note that in [13] a family
$$
    v^t =  |z_1|^2 + |z_2|^2, \ t\in \mathbb{R} \setminus \{0, 1, 2\}.
\eqno (5.12)
$$
of affine homogeneous real hypersurfaces was constructed.

  The value of $ t = 2 $  that is prohibited by this formula, corresponds
to the surface (2.6).

{\bf Remark.} Affine transformation groups of all the surfaces
(5.12) and (2.6) are 6-dimensional and can be easily written. One
can prove that there are only two (up to affine transformations)
second order surfaces $ v = 2 x_1^2 + |z_2|^2 $ and $ v^2 = |z_1|^2
+ |z_2|^2 $ that are affinely homogeneous manifolds of the type (0,
1/2) and have "reach" (more than 5-dimensional) affine
transformation groups.

\

\centerline{\textbf{5.5 Outline of the proof of the Theorem 4.5}}

\

  First, consider the case
$
   m_3 + i m_4 = 0.
$

  The similarity with the matrix (5.1) and
$ \lambda = m_1/2 $ reduces the bases of discussed algebras to the
form (5.8) where $ r = t_1 + m_1/2$ and conditions $ t_7 = t_8 = 0 $
are added.

  From the discussions above, the quadric
$
  v = 2 x_1^2 + |z_2|^2
$ appears when $ r = 0 $ and if $ r \ne 0 $ we get the logarithmic
surface $
  v = -\ln(1+x_1) + |z_2|^2
$.

  In the main case
$$
   m_3 + i m_4 \ne 0
\eqno (5.13)
$$
the proof of the Theorem 4.5 is based on the two following
statements.

{\bf Proposition  4.1} {\it If $ m_3 + i m_4 \ne 0 $ then any of the
algebras (3.5) is similar to the algebra with the basis
$$
E_1 =
 \left( \begin {array}{cccc} t_{{1}}&0&-i ( t_{{1}}m_{{1}}+2
\,{m_{{3}}}^{2}+2\,{m_{{4}}}^{2})/2 &0\\\noalign{\medskip}0&0&0&0
\\\noalign{\medskip}4\,i&0&2\,m_{{1}}&0\\\noalign{\medskip}0&0&0&0
\end {array}
\right) ,\eqno (5.14) $$ $$E_2 =
 \left( \begin {array}{cccc} 0&0&0&i\\\noalign{\medskip}0&0&0&0
\\\noalign{\medskip}0&0&0&0\\\noalign{\medskip}0&0&0&0\end {array}
 \right)
, E_3 =
 \left( \begin {array}{cccc} 2&0&0&0\\\noalign{\medskip}0&1&0&0
\\\noalign{\medskip}0&0&2&0\\\noalign{\medskip}0&0&0&0\end {array}
 \right)
,
$$
$$E_4 =
 \left( \begin {array}{cccc} 0&0&0&0\\\noalign{\medskip}0&i&0&0
\\\noalign{\medskip}0&0&0&0\\\noalign{\medskip}0&0&0&0\end {array}
 \right)
, \ E_5 =
 \left( \begin {array}{cccc} 0&0&0&0\\\noalign{\medskip}0&0&0&0
\\\noalign{\medskip}0&0&0&1\\\noalign{\medskip}0&0&0&0\end {array}
 \right)
.
$$
}

\

{\bf Proposition 4.2.} {\it
  Integral variety of any algebra (5.14)
is affine equivalent to a surface of the following family:

$$
   v = x_1^{(1 - \alpha)} |z_2|^{2\alpha}, \
        \alpha \in \mathbb{R} \setminus \{0,1\}.
$$
}

  We give a brief comment to the proof of Proposition 4.2.

   Due to the simple form of the
matrices $ E_2, E_4, E_5 $, the study of the system of five
equations corresponding to basic fields of algebra (5.14) becomes
significantly simpler. We need to solve just two equations
$$
  (t_1 x_1 + A G)\frac{\partial G}{\partial x_1}  = 4 x_1 + 2 m_1 G,
\quad
    x_1 \frac{\partial G}{\partial x_1} +
      s \frac{\partial F}{\partial  s } =  G
\eqno (5.15)
$$
for an analytic function $ G(x_1,s) $. Here
$$
   A = \frac12( t_1 m_1 + 2 \ m_3^{2}+ 2\, m_4^{2}),
\quad s = |z_2|^2,
$$
and the desired equation of the surfaces has the form $ v = G(x_1,
|z_2|^2) $.

   The common solution of the second equation of (5.15) can be presented
(outside the set $ x_1 = 0 $) in the form
$$
    G = x_1 H(s/x_1)
\eqno (5.16)
$$
with an arbitrary analytic function $ H $. Then the first equation
in (5.15) takes the form of ODE
$$
   \xi (t_1 + A H) H' = A H^2 + B H - 4,
\eqno (5.17)
$$
where
$$
   \xi = \frac{s}{x_1}, \quad B = (t_1 - 2 m_1).
$$

  The following consideration of two cases corresponding to the equality or
inequality to zero of a parameter $ A $, leads to a common formula
(2.9) for the unknown homogeneous surfaces but with different
domains of variation of the parameter $ \alpha $.

   Note that the formulas which connect the parameters of the homogeneous surfaces
families have extremely cumbersome form in different representations
of these families. This fact, known for a long time (see [3]), also
appears in the proof of the Theorem 4.5.

\

\centerline{\textbf{5.6 Outline of the proof of the Theorem 4.6}}

\

  Here it is also convenient, as in the proof of the Theorem 4.5, to consider
two cases related to the parameter
$
  m_3 + i m_4.
$
   In the first simple case
$
  m_3 = m_4 = 0
$ the standard similarity (5.1) with $ \lambda = m_1/2 $ transforms
the basis (3.6) to the basis of 5-dimensional subalgebra of the
algebra (2.16). Consequently, in this case, the desired homogeneous
surface is the standard quadric (2.15).

   In the general case
$
  m_3 + i m_4 \ne 0,
$ we note that there is a natural symmetry between $m_3$ and $m_4$
in algebras of the family (3.6). It is sufficient therefore  to
discuss only a situation
$$
m_3 \ne 0. \eqno (5.18)
$$

  The consideration of the dual case $ m_4 \ne 0 $ gives the result which
is affine
equivalent to the one obtained with the assumption (5.18).

   For any of the algebras (3.6) the spectrum of the matrix
$$
E_3 = \left(%
\begin{array}{cccc}
  \frac{3}{2}m_3 & -m_1 & \frac{i}{4}m_1m_3 & 0 \\
  0 & m_3 & 0 & 1 \\
  0 & 2i & 2m_3 & 0 \\
  0 & 0 & 0 & 0 \\
\end{array}%
\right),
$$
is presented by the set
$
   \left(\frac32 m_3, m_3, 2m_3, 0 \right).
$

Taking into account the condition (5.18), we can diagonalize this
matrix. Due to the similarity with the matrix
$$
   S =
 \left( \begin {array}{cccc} 1&m_{{1}}&im_{{1}}&m_{{1}}
\\
\noalign{\medskip}0&m_{{3}}&0&2\,m_{{3}}\\
\noalign{\medskip}0&-2\,i&
2&-2\,i\\
\noalign{\medskip}0&0&0&-2\,{m_{{3}}}^{2}\end {array}
 \right)
$$
consisting of the eigenvectors of $ E_3 $, one can switch to a new
algebra.

  Using another similarity with diagonal matrix
$$
   D = diag \left(\frac{m_{3}^{2}+ m_{4}^{2}}{m_3}, \frac{m_3 +im_4}{m_3},
\frac{m_3^{2}+ m_4^{2}}{m_3^{2}},1 \right)
$$
any of the algebras under consideration becomes the basis
independent of the parameters $ m_3, m_4 $
$$
E_1=
 \left( \begin {array}{cccc} 0& 1
&0&0\\
\noalign{\medskip}0&0&0&0\\
\noalign{\medskip}2\,i&0&0&0
\\
\noalign{\medskip}0&0&0&0\end {array} \right), \ E_2 =
 \left( \begin {array}{cccc} 0&0&0&i\\
\noalign{\medskip}0&0&0&0
\\
\noalign{\medskip}0&0&0&0\\
\noalign{\medskip}0&0&0&0\end {array}
 \right)
, \eqno (5.19)$$ $$E_3 =
 \left( \begin {array}{cccc} 3&0&0&0\\
 \noalign{\medskip}0&2&0&0
\\
\noalign{\medskip}0&0&4&0\\
\noalign{\medskip}0&0&0&0\end {array}
 \right)
, E_4 = \left( \begin {array}{cccc} 0&0&0&0\\
\noalign{\medskip}
0&0&0& -2\,i \\
\noalign{\medskip}0&1&0&0
\\
\noalign{\medskip}0&0&0&0\end {array} \right) , \ E_5 =
 \left( \begin {array}{cccc} 0&0&0&0\\
 \noalign{\medskip}0&0&0&0
\\
\noalign{\medskip}0&0&0&1\\
\noalign{\medskip}0&0&0&0\end {array}
 \right)
.
$$

   The system of partial differential equations corresponding to such algebra
contains three equations (taking into account the trivial form of
the matrices $ E_2 $ and $ E_5 $)
$$
     x_2 \frac{\partial F}{\partial x_1} =  2 x_1,
\qquad
     2 \frac{\partial F}{\partial y_2} = - y_2,
\eqno (5.20)
$$
$$
   3 x_1 \frac{\partial F}{\partial x_1} +
   2 x_2 \frac{\partial F}{\partial x_2} +
   2 y_2 \frac{\partial F}{\partial y_2} = 4 F
$$
with unknown function $ F(x_1, x_2, y_2 ) $.

  After step-by-step integration of these equations we obtain the
equation of a homogeneous surface in the form
$$
    F = \frac{x_1^2}{x_2} - \frac14 y_2^2  +  C x_2^2,
\eqno(5.21)
$$
where $ C $ is an arbitrary constant.

   Tracing the movement of the origin (that belongs to the unknown
homogeneous surface) under the intermediate transformations, we can
determine the value of the constant $ C = -1 / 4 $ in the formula
(5.21).

   As a result, we arrive at equation (2.7) defining
(up to affine transformation) a homogeneous surface in this case.

{\bf Remark.} One can easily show that for $ |C| \ne 1/4 $ the
surface (5.21) is affine homogeneous manifold of the type $ (1/2,
(1+ 4C)/(1-4C) $. All such (affinely different!) surfaces are the
integral manifolds of the same algebra with basis (5.19). Note also
that they are holomorphically equivalent to each other.

\

\centerline{\textbf{5.7 Outline of the proof of the Theorem 4.7}}

\

   As in the previous discussion, we choose one complex
parameter $ t_3 + i t_4 $ from the real parameter quadruplets $
   m_1, t_1 $,
$ t_3, t_4, $ on which the family (3.7) depends.

   Note that the basis of the form (3.7) for $ t_3 + i t_4 = 0 $ and the
basis of the form (3.5) considered above identically coincide when $
m_3 + i m_4 = 0 $.
     Consequently, the homogeneous surfaces obtained in these two subcases,
i.e. $
  v = 2 x_1^2 + |z_2|^2
$ and $
  v = -\ln(1+x_1) + |z_2|^2
$ are the same in both cases.

\

{\bf Proposition 4.3.} {\it If $ t_3^2 + t_4^2 \ne 0 $ then the
basis of any algebra of the discussed family can be transformed by
matrix similaritiy to the form
$$
E_1 =
\left( \begin {array}{cccc} t_{{1}}&-iA&0&0\\
\noalign{\medskip}4\,i&2\,m_{{1}}&0&0\\\noalign{\medskip}0&0&0
&0\\\noalign{\medskip}0&0&0&0\end {array} \right) , E_2 = \left(
\begin {array}{cccc} 0&0&1&0\\\noalign{\medskip}0&0&0&0
\\\noalign{\medskip}0&0&0&0\\\noalign{\medskip}0&0&0&0\end {array}
 \right)
, \eqno (5.22)
$$
$$E_3 =
 \left( \begin {array}{cccc} 0&0&0&0\\\noalign{\medskip}0&0&1&0
\\\noalign{\medskip}0&0&0&0\\\noalign{\medskip}0&0&0&0\end {array}
 \right)
, E_4 =
 \left( \begin {array}{cccc} 1&0&0&0\\\noalign{\medskip}0&1&0&0
\\\noalign{\medskip}0&0&-1&0\\\noalign{\medskip}0&0&0&0\end {array}
 \right)
, E_5 =
 \left( \begin {array}{cccc} i&0&0&0\\\noalign{\medskip}0&i&0&0
\\\noalign{\medskip}0&0&i&0\\\noalign{\medskip}0&0&0&0\end {array}
 \right),
$$
}

where $A = t_3^2 + t_4^2 + \frac{t_1m_1}2$.

  Outline of the proof of Proposition 4.3 is similar to the one from
Proposition 4.1.
   In this case we have the symmetric dependence of the bases of the algebras
on parameters $ t_3 $ and $ t_4 $ (similar to the parameters $ m_3 $
and $ m_4 $ in section 4.5). So, in fact it is sufficient to
consider only the subcase $ t_4 \ne 0 $.

    In this subcase at first we
 diagonalize (by the similarity) a matrix $ E_4 $, and then
consider linear combinations of obtained basis matrices. Such
actions allow us to "improve" \ the basis in total.

    At the next step of the proof of the Theorem 4.7 we simplify the upper left $ (2 \times 2) $-block $ e_1 $ of the matrix $ E_1 $ from the
basis (5.22). If $
   A = t_3^2 + t_4^2 + \frac{t_1m_1}2
$
 then
$
   e_1 =
\left( \begin {array}{cc}
t_1 & -iA \\
4i  & 2 m_1
\end {array}
 \right)
\quad $ has two real eigenvalues of opposite signs
$$
  \lambda_{1,2} = \frac12 ((t_1 + 2m_1) \pm B), \
\mbox{ \ where } \
  B = \sqrt{(t_1 + 2m_1)^2 + 16(t_3^2 + t_4^2) }.
\eqno (5.23)
$$

  Consequently, the matrix $ E_1 $ is diagonalizable
and after an appropriate matrix similarity transformation it has the
form $$
  E_1^* =
 \left( \begin {array}{cccc} 1 &0&0&0\\
0& \alpha &0&0\\
0&0&0&0\\
0&0&0&0\end {array}
 \right), \quad \alpha =  \frac{\lambda_1}{\lambda_{2}} < 0.
$$

We note here two important facts:

1) based on the formula (5.23), it is easy to show that the ratio $
\alpha $ takes all the negative values depending on the parameters $
t_1, m_1, t_3, t_4 $;

2) it can be also shown that the similarity transformation
diagonalizing $ E_1 $ does not actually change the rest of the basis
matrices of the discussed algebras. \

   Now it remains to integrate algebras with the bases
$$
  E_1 =
 \left( \begin {array}{cccc} 1 &0&0&0\\
0& \alpha &0&0\\
0&0&0&0\\
0&0&0&0\end {array}
 \right)
, \ E_2 = \left( \begin {array}{cccc}
0&0&1&0\\\noalign{\medskip}0&0&0&0
\\\noalign{\medskip}0&0&0&0\\\noalign{\medskip}0&0&0&0\end {array}
 \right)
,\eqno (5.24)$$ $$E_3 =
 \left( \begin {array}{cccc} 0&0&0&0\\\noalign{\medskip}0&0&1&0
\\\noalign{\medskip}0&0&0&0\\\noalign{\medskip}0&0&0&0\end {array}
 \right)
, E_4 =
 \left( \begin {array}{cccc} 1&0&0&0\\\noalign{\medskip}0&1&0&0
\\\noalign{\medskip}0&0&-1&0\\\noalign{\medskip}0&0&0&0\end {array}
 \right)
, E_5 =
 \left( \begin {array}{cccc} i&0&0&0\\\noalign{\medskip}0&i&0&0
\\\noalign{\medskip}0&0&i&0\\\noalign{\medskip}0&0&0&0\end {array}
 \right)
.
$$

\

{\bf Proposition 4.4.} {\it
 Integral variety of any algebra with the basis (5.24)
is affinely equivalent to one of the surfaces of the following
family:
$$
    Re  (\bar z_1 w) = C (Re (z_1 \bar z_2))^{\alpha}, \quad
         \alpha \in (- \infty,0), \ \quad C \in \mathbb{R}.
\eqno (5.25)
$$
}

   To prove the Proposition 4.4, at first we transform each basis
   matrix (5.24) into a lower triangular matrix that is convenient
   for the following integration. It can be easily done by the change of
variables
$$
    z_1 = z_2^*, \quad z_2 = w^*, \quad w= z_1^*
\eqno (5.26)
$$
(or, what is the same, by another similarity transformation).

   After this transformation, the system of partial differential equations corresponding
   to the modified basic fields and to the desired equation
$ v = F (z, \bar z, u) $ of homogeneous surface, takes the form
$$
E1: \  x_2 \frac{\partial F}{\partial x_2} +
  y_2 \frac{\partial F}{\partial y_2} +
  \alpha (u \frac{\partial F}{\partial u} - F) = 0,
$$
$$
E_2: \ x_1 \frac{\partial F}{\partial x_2} +
  y_1 \frac{\partial F}{\partial y_2}  = 0,
\quad E_3: \  - y_1 \frac{\partial F}{\partial u} = x_1, \eqno
(5.27)
$$
$$
E_4: \  - (x_1 \frac{\partial F}{\partial x_1} +
           y_1 \frac{\partial F}{\partial y_1}) +
          (x_2 \frac{\partial F}{\partial x_2} +
           y_2 \frac{\partial F}{\partial x_2}) -
          (u \frac{\partial F}{\partial u} - F) = 0
$$
$$
E_5: \  (- y_1 \frac{\partial F}{\partial x_1} +
           x_1 \frac{\partial F}{\partial y_1} ) +
        (- y_2 \frac{\partial F}{\partial x_2} +
           x_2 \frac{\partial F}{\partial y_2} ) -
         F \frac{\partial F}{\partial u}  = u.
$$

  Solving the equations corresponding to the fields
$ E_3 $, $ E_2 $ and $ E_4 $ sequentially, we obtain the formula
$$
  F = - \frac{x_1 u }{y_1} + \frac1{x_1}\cdot
      \varphi\left(\frac{y_1}{x_1}, x_1 y_2 - x_2 y_1 \right),
$$
where $ \varphi $ is an arbitrary analytic function of two
arguments.

   Denoting
$$
    \xi = \frac{y_1}{x_1}, \quad r = x_1 y_2 - x_2 y_1,
$$
we obtain the remaining two equations of (7.10) in the form
$$
   r \frac{\partial \varphi}{\partial r} = \alpha \varphi,
   \quad
   (\xi^2 + 1)\cdot (\frac{\partial \varphi}{\partial \xi} + \varphi) = 0.
\eqno (5.28)
$$

   The solving of this system leads us to the
equation of the required homogeneous surface of the form
$$
    v + \frac{x_1 u }{y_1} = C \frac1{y_1} (x_1 y_2 - x_2 y_1)^{\alpha}
$$
with an arbitrary real constant $ C $.

  In the complex variables, this equation can be rewritten as
$$
   Re\, (\bar z_1 w) = C \left( Im\, (\bar z_1 z_2)\right)^{\alpha}.
$$

The change of variable $ z_2 =-i z_2 ^* $ let us use the symbol $Re$
instead of the symbol $Im$ in the right hand side of the equation.
Proposition 4.4 is proved.

  Note now that for zero constant $ C $ the last equation
defines a degenerate surface in Levi sense. If $C$ is a nonzero
constant, we can assume that $C = 1$ due to the scaling of variable
$ w $. Thus, we prove the Theorem 4.7.

{\bf Remark.} It is easy to verify directly that the equation (2.10)
defines an affine homogeneous surfaces for all values of $ \alpha $.
However, for $ \alpha = 0 $ the surface (2.10) degenerates in the
sense of Levi. Additional simple analysis of power series shows that
for
  $ \alpha> 0 $
(2.10) describes the surfaces with an indefinite Levi form, that also
do not satisfy the strictly pseudo-convexity condition.

\

\

\centerline{{\bf References}}

\

\begin{enumerate}
    \item Cartan E. Sur la geometrie pseudoconforme des hypersurfaces
         de deux va\-riables complexes // Ann. Math. Pura Appl.,
         (4) 11 (1932), P. 17 - 90 (Oeuvres II, 2, 1231 - 1304 ).
    \item Azad H., Huckleberry A., Richthofer W. Homogeneous CR manifolds //
        J. Reine und Angew. Math. Bd. 358 (1985), P. 125 - 154.
    \item Eastwood M., Ezhov V.V. On affine normal forms and a classification
     of homogeneous surfaces in affine three-space // Geom Dedicata, 1999,
     V. 77, P. 11 - 69.
    \item Ezhov, V, Isaev, A., Schmalz, G., Invariants of elliptic and
   hyperbolic CR-structures of codimension 2, Int. Journ. of Math.,
   V. 10 (1999), N 1, pp. 1-52.
    \item Isaev, A., Analogues of Rossi's map and E. Cartan's
   classification of homogeneous strongly pseudoconvex 3-dimensional
   hypersurfaces, Journal of Lie Theory, vol. No. 16 (2006), pp. 407-426.
    \item Fels, G.  Classification of Levi degenerate homogeneous
      CR-manifolds in dimension 5 / G. Fels, W. Kaup // Acta Math. -- V.
      210(2008). -- P. 1 - 82.
    \item Merker, J., Lie symmetries and CR Geometry,
   J. of Math. Sciences, 154 (6), pages 817--922, 2008.
    \item Lamel, B., Mir N., Zaitsev D., Lie group structures on automorphism groups
   of real-analytic CR manifolds, Amer. Journ. of Math., V. 130 (2008), No. 6,
   pp. 1709-1726.
    \item Beloshapka, V.K., Kossovskiy I.G., Homogeneous hypersurfaces in
   $\mathbb{C}^{3}$, associated with a model CR-cubic, J. Geom. Anal.
   V. 20 (2010), No. 3, P. 538 - 564.
    \item Beloshapka, V.K., Kossovskiy I.G., Classification of homogeneous CR-
   manifolds in dimension 4, J. Math. Anal. Appl.,  V. 374 (2011),
   No. 2. P. 655-672.
    \item Loboda, A.V., Affinely Homogeneous Real Hypersurfaces of $ \mathbb{C}^2 $.
   Funkts. Anal. Prilozh., 47:2 (2013), 38-54.
    \item Loboda, A. V., "Homogeneous Strictly Pseudoconvex Hypersurfaces in C3 with
   Two-Dimensional Isotropy  Groups," Mat. Sb. 192 (12), 3-24 (2001)
   [Sb. Math. 192, 1741-1761 (2001)].
    \item Loboda, A. V. and Khodarev, A. S., "On a Family of Affine-Homogeneous Real
   Hypersurfaces of a Three-Dimensional Complex Space," Izv. Vyssh. Uchebn.
   Zaved., Mat., No. 10, 38-50 (2003) [Russ. Math. 47 (10),
   35-47 (2003)].
    \item A. V. Loboda and T. T. D. Nguyen, On the Affine Homogeneity of Tubular
   Type Surfaces in $ \mathbb{C}^3 $, Tr. Mat. Inst. im. V.A. Steklova,
   Ross. Akad. Nauk 253, 102-119 (2012)
   [Proc. Steklov Inst. Math. 279, 93-110 (2012)]
    \item A. V. Loboda, On complete description of affine-homogeneous real
   hypersurfaces of tubular type in $\mathbb{C}^3 $. Abstracts of Voronezh
   Winter Mathematical School, (Voronezh, 2013), 144-145.
    \item T. T. D. Nguen, Affine-Homogeneous Real Hypersurfaces of Tubular Type in
   $\mathbb{C}^3 $. Mat. Zametki, 94:2 (2013),  246-265.
    \item Shabat, B.V., Introduction to Complex Analysis, Part 2:
    Functions of Several Variables (Nauka, Moscow, 1976) [in Russian].
    \item Mubarakzyanov, G. M., Classification of real structures of Lie algebras
    of fifth order. Izv. Vyssh. Uchebn. Zaved. Mat., 1963, no. 3 (34), 99-106.
    \item Sabzevari M., Hashemi A., Alizadeh B.M., Merker J. Applications
    of differential algebra for computing Lie algebras of
    infinitesimal CR-automorphisms, arXiv:1212.3070 (December 2012)
\end{enumerate}

\end{document}